\documentclass[11pt]{amsart}
\usepackage{amssymb,amsfonts,amsthm,amsmath}
\usepackage[english]{babel}
\usepackage{mathrsfs}

\newtheorem{theorem}{Theorem}[section]
\newtheorem{lemma}[theorem]{Lemma}
\newtheorem{corollary}[theorem]{Corollary}
\newtheorem{definition}[theorem]{Definition}
\newtheorem{proposition}[theorem]{Proposition}
\newtheorem{remark}[theorem]{Remark}

\newcommand{\field}[1]{\mathbb{#1}}
\newcommand{\Cal}[1]{\mathcal{#1}}
\newcommand{\Scr}[1]{\mathscr{#1}}
\newcommand{\la}{{\langle}}
\newcommand{\ra}{{\rangle}}
\newcommand{\absgamma}{{|\gamma|}}
\newcommand{\ba}{{\bf a}}
\newcommand{\bB}{{\bf B}}
\newcommand{\bg}{{\bf g}}

\newcommand{\bo}{{\bf 0}}
\newcommand{\bsym}{{\bf Sym}}

\newcommand{\clos}{{\rm clos}}

\newcommand{\Q}{\field{Q}}
\newcommand{\R}{\field{R}}
\newcommand{\crit}{{\rm crit}}
\newcommand{\dr}{\partial_r}

\newcommand{\drf}{\partial_r f}

\newcommand{\Hs}{{\rm Hess}}
\newcommand{\Hsg}{{\rm Hess}^{\bg}}
\newcommand{\nf}{{\nabla \! f}}
\newcommand{\ngf}{{\nabla_{\! \bg} f}}
\newcommand{\ngu}{{\nabla_{\! \bg} u}}
\newcommand{\np}{{\nabla^{'}}}
\newcommand{\Rn}{\R^n}
\newcommand{\Sr}{{\bf S}}
\newcommand{\Srn}{\Sr^{n-1}}

\newcommand{\Symn}{{\rm Sym (\Rn)}}
\newcommand{\nuf}{{\nu_f}}
\newcommand{\nugf}{{\nu^\bg_f}}
\newcommand{\ud}{{\rm d}}
\newcommand{\ve}{{\varepsilon}}

\newcommand{\bs}{ {\tiny $\blacksquare$} \\}

\numberwithin{equation}{section}

\title[]{On Hessian limit directions along non-oscillating gradient trajectories}
\author{Vincent Grandjean}
\address{{\it Temporary Address:} V. Grandjean, Fields Institute, 222 College Street,
Toronto, Ontario,  M5T 3J1, Canada}
\email{cssvg@bath.ac.uk}
\thanks{}
\date{\today}
\begin{document}
\maketitle
\begin{abstract}
Given a non-oscillating gradient trajectory $\absgamma$ of a real analytic function $f$,
we show that the limit $\nu$ of the secants at the limit point $O$ of $\absgamma$ along the trajectory
$\absgamma$ is an eigen-vector of the limit of the direction of the Hessian matrix $\Hs (f)$ at $O$
along $\absgamma$. The same holds true at infinity if the function is globally subanalytic. We also deduce
some interesting estimates along the trajectory. Away from the ends of the ambient space, this property is
of metric nature and still holds in a general Riemannian analytic setting.
\end{abstract}

\bigskip
\section{introduction}
Let $f :(\Rn,O) \to (\R,0)$ be the germ at $O$ of a real analytic function singular at $O$.
The behavior of the germ at $O$ of gradient trajectories of the function $f$ with $\omega$-limit
point $O$ is barely understood. Thom's Gradient conjecture, proved by
Kurdyka, Mostowski and Parusi\'nski \cite{KMP}, states that at $O$ the trajectory has a tangent. In whole
generality nothing else is known. The plane smooth case is very well understood, and (bounded) plane
gradients curves are definable in the Pfaffian closure of $\R_{an}$. But a few special cases trying to address
the non-spiraling/non-oscillating behavior of such trajectories \cite{Sa,FS,GS}, and the main result of
\cite{Mo}, that is all that is known.

\medskip
The aim of this note is to report an intriguing and not expected result,
along gradient trajectories that are non-oscillating at their limit point. We recall that a gradient 
trajectory
$\gamma$ is (\em semi-analytically\em) \em non-oscillating \em at its limit point $O\in \crit (f)$ if, 
given any semi-analytic set $H$ (germ at $0$), either
$\gamma \subset H$ or $\gamma \cap H = \emptyset$ (as germ at $O$). Our main result is

\medskip\noindent
{\bf Theorem \ref{Theorem:main}.} \em
Let $\gamma$ be a gradient trajectory of a real analytic function $f: (\R^n,O) \to (\R,0)$.
Assume that $\gamma$ is non-oscillating at its $\omega$-limit point $O\in \crit (f)$.
Let $\nu\in\Sr^{n-1}$ be the limit at $O$ of the secant along $\gamma$. \\
1) The oriented direction of the matrices $\Hs (f)$ has a (non-zero) limit $\Cal{H}$ at $O$ along $\gamma$. 

\smallskip\noindent
2) $\displaystyle{\lim_{\absgamma\ni x \to O}\frac{\Hs (f)(x)\cdot\nf (x)}{|\Hs (f)(x)\cdot\nf (x)|}} = -\nu$.
Thus the unit vector $\nu$ is an

\smallskip\noindent 
eigen-direction of $\Cal{H}$.\em

\medskip\noindent
This result is of metric nature, and although it is here just stated for the Euclidean metric it still holds
in any analytic Riemannian setting (see Section \ref{Section:riemannian}).
This result is somehow providing a clue along non-oscillating gradient trajectories about
where a gradient trajectory is ending-up. Note that this result was not even known in the plane case,
where it applies.

\medskip
The paper is organized as follows.

Section \ref{Section:background} introduces the few useful definitions, notions and notations we will deal
with.

In Section \ref{Section:setting} we present quickly the framework which we will operate in, and recall
\cite{KMP} an interesting technical result for us.

Section \ref{Section:hessian-and-gradient} is aimed at giving motivations to this paper, in stating and
re-proving results about links between gradient trajectories and Hessian matrices that should be 
well known folklore. These special cases suggest two possible different
interpretations of these connections between Hessian matrices and gradient trajectories.

The main result is stated and proved in Section \ref{Section:main-result}.

It is followed by a few corollaries in Section \ref{Section:AFUC} presenting some estimates along the
trajectory.

Section \ref{Section:infinity} deals with the analog result at the limit point at infinity of a
non-oscillating (at infinity) unbounded gradient trajectory for of globally subanalytic real analytic 
function.
We also provide the corresponding estimates expected from those found in the previous section.

Section \ref{Section:riemannian} we state and prove our main result in the general (real analytic)
Riemannian case.

\medskip\noindent
{\bf Thanks.}  I am very grateful to the Fields Institute for its support and the working
conditions provided since this work started during the Fields Institute's Thematic Program on \em
O-minimal Structures and Real Analytic Geometry January-June 2009. \em I thank F. Sanz a lot
for having suggested to jump from the plane case to the non-oscillating case. I owe many thanks
to K. Kurdyka for precious comments, suggestions and encouragements.
%
%
%
%
%
%
%
%
%
%
%
%
%
%
%
%
%
%
\section{Background - Notations}\label{Section:background}
Let $\gamma : I \to \Rn$ be $C^1$ parameterized curve over a connected interval $I$ with
non-empty interior.
We will write $\absgamma$ for the support of the parameterized curve, that is
\begin{center}
\vspace{4pt}
$\absgamma = \{\gamma (t), t\in I\}$.
\vspace{4pt}
\end{center}

\begin{definition}\label{Definition:oscillating}
The curve $\absgamma$ is non-oscillating if for any semi-analytic subset $H$, the intersection
$H\cap\absgamma$ has finitely many connected components.
\end{definition}

Since we are interested in a local problem we prefer using the following local version:
\begin{definition}\label{Definition:oscillating-at-O}
Assume that the origin $O$ belongs to $\clos(\absgamma)\setminus\absgamma$.
The curve $\absgamma$ is non-oscillating at $O$, if for any semi-analytic subset $H$, the germ at
$O$ of the intersection $H\cap\absgamma$ is either empty or is the germ of $\absgamma$ at $O$.
\end{definition}

These notions were explored for real analytic vector fields \cite{Sa,CMS,FS} on $3$-manifolds to study
their asymptotic dynamics at their singular points in certain situations, namely to discuss spiraling.

\bigskip
The pointed vector space $\Symn^*$ of non-zero real symmetric matrices of $\Rn$ is an open cone
on the unit sphere of $\Symn$ with vertex the null matrix, on which the group $(\R_{>0},\cdot)$ acts
naturally, smoothly and semi-algebraically by the homotheties of positive ratio.
The resulting quotient space $\bsym (n)$ is the space of {\it oriented directions of the real symmetric
$n$-matrices} and is a real algebraic manifold isomorphic to $\Sr^{\frac{n(n+1)}{2}-1}$.
The quotient mapping $\sigma : \Symn^* \to \bsym (n)$  is a smooth semi-algebraic submersion.

\smallskip
Let $f$ be a $C^2$ function defined on some open set $U$ of $\Rn$ equipped with some coordinates.
We define,  $\Hs (f) =[h_{i,j}]$ to be the Hessian matrix of $f$, and the
function $h$ on $U$ as $h(x) = [\sum_{i,j}h_{i,j}(x)^2]^{1/2}$. We will write $H(f)$ for the matrix
$h^{-1}\Hs (f)$ where it is defined. We obviously get $\sigma(\Hs (f)) = \sigma (H(f))$.
Note that $H(f)$ is not $\sigma (H(f))$, but has a similar property: let $(\Hs (f)(x_k))_k \in \Symn^*$
be a sequence converging to the null matrix. The accumulation values of the sequence
$(H(f)(x_k))_k$ are bounded and will not all be the null matrix.
In particular if the sequence $(H(f)(x_k))_k$ has a
unique accumulation value $\Cal{H}$, the sequence $\sigma (H(f)(x_k))_k$ then converges to
$\sigma (\Cal{H}) = \lambda \Cal{H}$ for some positive real number $\lambda$. With an obvious but small 
abuse of language we will also call such a limit $\Cal{H}$ an \em oriented limit direction \em of the
sequence $(H(f)(x_k))_k$.

\bigskip
Last a short word on notations.

\smallskip
Unless explicitly mentioned $\Rn$ comes equipped with the Euclidean structure and we will write
$\la -,-\ra$ for the scalar product and $|-|$ for the associated norm.

%
%
%
%
%
%
%
%
%
%
%
%
%
%
%
%
%
%
\section{Setting}\label{Section:setting}
Let $f:U \to \R$ be a real analytic function defined on an open subset $U \subset \Rn$ containing the
origin $O$, such that $O \in \crit (f)$ the critical locus of $f$. We also assume $f(O) = 0$.

Let $\nuf$ be the unitary gradient vector field defined on $\Rn \setminus \crit (f)$.
But its orientation, this vector field does depend only on the foliation by the connected
components of the levels of the function, not on the chosen function.

Let us consider the following gradient differential equation:
\begin{center}
\vspace{4pt}
$\dot{\gamma}(s) = \nuf (\gamma (s))$ with $\gamma (0) = x_0$ and $f(x_0) < 0$.
\vspace{4pt}
\end{center}
Given any maximal non stationary solution $\gamma$, let $\absgamma$ be the underlying subset
of the parameterized integral curve, namely
\begin{center}
\vspace{4pt}
$\absgamma =\{ \gamma (s), s\geq 0 \}$.
\vspace{4pt}
\end{center}

\begin{definition}
We call $\absgamma$ a half-trajectory of the gradient field associated with the function $f$.
\end{definition}
For short we may say trajectory instead of half-trajectory.

\medskip
Assume we are given a half-trajectory $\absgamma$, we assume that the function $ s\to f(\gamma (s))$
increases to $0$ and that the omega limit set $\omega (\absgamma)$ contains $O$. \\
Let $l(\absgamma)$ be the length of $\absgamma$ between $\gamma (0)$ and the origin $O$. We know since
\L ojasiewicz that it is finite \cite{Loj}, and so $\omega (\absgamma) =\{O\}$.\\
The next result is the first systematic result on the behavior of gradient trajectories
at their limit point, and so far the only such one within the category of analytic functions.
\begin{theorem} [\cite{KMP}]\label{Theorem:gradient-conjecture}
{\it The length of the radial projection  $t\to \frac{\gamma(t)}{|\gamma (t)|}$ onto $\Sr^{n-1}$
of the gradient curve $\gamma$  is finite.}
\end{theorem}
\noindent
The straightforward and looked for consequence of Kurdyka, Mostowski \& Parusi\'nski's
result is that the well-known \em Thom's Gradient Conjecture \em holds true:

\smallskip\noindent
{\bf Thom's Gradient Conjecture.} \em There exists $\nu \in \Srn$ such that\em
\begin{center}
\vspace{4pt}
$\displaystyle{\lim_{s \to l(\absgamma)}\frac{\gamma (s)}{|\gamma (s)|} = \nu}$ or equivalently
$\displaystyle{\lim_{\absgamma\ni x \to O}\frac{x}{|x|} =  \nu}$.
\vspace{4pt}
\end{center}

\medskip
Writing $\dr$ for the radial vector field,
the gradient vector field $\nf$ then decomposes into two orthogonal components  as $\nf = \dr f \dr +
\np \! \! f$ where $\dr f \dr$ is its radial component and  $\np \!\! f$ is its spherical part.

To prove their result, Kurdyka \& Als used fine estimates along the trajectory. In particular we mention
the following result that will prove to be useful for our purpose:
\begin{proposition}[\cite{KMP}]\label{Proposition:exponent-vca}
1) There exists a rational number $m := m(\absgamma)$
and a positive real number $\ba := \ba(\absgamma)$ such that along $\absgamma$
\begin{center}
\vspace{4pt}
$ f = - \ba r^m + o(r^m)$ and $\dr f =  - m\ba r^{m-1} + o(r^{m-1})$
\vspace{4pt}
\end{center}
2) There exists a finite subset $S\subset\Q_{>0}\times \R_{>0}$ of pairs $(l,b)$, 
such that for any trajectory $\alpha$, such that $\omega (|\alpha|) =\{O\}$
and $f$ is negative along $|\alpha|$, the pair $(m(|\alpha|),\ba(|\alpha|)$ 
as in point 1) corresponding to $|\alpha|$ belongs to $S$.
\end{proposition}
%
%
%
%
%
%
%
%
%
%
%
%
%
%
%
%
%
%
%
%
%
%
%
%
%
%
%
%
%
\section{Hessian and Gradient: special cases}\label{Section:hessian-and-gradient}
This section is devoted to present some examples exhibiting strong connections between
limits of direction of Hessian matrices and gradient vector fields. This material
should be well known and at least part of the folklore of this subject.
There are somehow similar suggestions in \cite{Mo}.

\medskip
1) We decompose the real analytic function $f$ as the sum of its homogeneous parts at $O$, namely
\begin{center}
\vspace{4pt}
$f= f_p + f_{p+1}+ \ldots$
\vspace{4pt}
\end{center}
Let $\Hs (f)$ be the Hessian matrix of $f$.

We recall that for each $k\geq p$ the following formulae hold
\begin{center}
\vspace{4pt}
$k(k-1)f_k = (k-1) \la \nf_k, r\dr\ra =  \la \Hs (f_k)\cdot r\dr, r\dr\ra $,
\vspace{4pt}
\end{center}
Let $\beta:\Sr^{n-1} \times \R_+ \to \Rn$ be the spherical blowing-up with center $O$ defined as 
$(u,r) \to ru$.
Up to a multiplication by $p(p-1)r^{2-p}$ the gradient differential equation now reads:
\begin{center}
\vspace{4pt}
\begin{tabular}{rcl}
\vspace{4pt}
$\dot{r} $ & = & $r\la \Hs (f_p)(u)\cdot u,u\ra +O(r^2)$ \\
$\dot{u} $ & = & $\Hs (f_p)(u)\cdot u - \la \Hs (f_p)(u)\cdot u,u\ra u +O(r).$
\vspace{4pt}
\end{tabular}
\end{center}
The vector field corresponding to this differential equation is called
the {\it divided gradient vector field}, and is denoted $\xi_f$. It extends analytically on $\Sr^{n-1}
\times \R$.

Thus if $v\in \Srn$ is not an eigen-vector of $\Hs(f_p)(v)$, then $\dot{u}|_{r=0,u=v} \neq 0$.
Thus any limit point $\nu$ of a gradient trajectory is such that $\nu$ is an eigen-vector of $\Hs(f_p)(\nu)$.
(Which is not such an interesting information once $\Hs(f_p)(\nu)$ is the null-matrix.)

\medskip
The following elementary result should be known, although we never came across it this way.
\begin{proposition}\label{Proposition:hessian-eigen-homogeneous}
Let $f$ be a homogeneous polynomial of degree $p\geq 2$. \\
Then along any gradient trajectory $\absgamma$ with $\omega$-limit set $O$ with limit of secants $\nu$,
we find $f(\nu) \neq 0$ so that $\nu$ is an eigen-vector of the non null-matrix $\Hs(f)(\nu)$ associated
with a non-zero eigen-value, and so the vector field $\xi_f$ is reduced at $\nu$, that is has a non-nilpotent
linear part.
\end{proposition}
\begin{proof}
Let $\absgamma$ be a trajectory of $\nf$ with $\omega$-limit point $O$ and limit of secants $\nu$.
Let $r(s) = |\gamma (s)|$ and let $\nu (s) =  \nu_f (\gamma (s))$, then $f \circ \gamma (s) = r(s)^p
f (\nu (s))$. Proposition \ref{Proposition:exponent-vca} gives along $\absgamma$,  $r\dr f  = m f +
o(r^m) = -m\ba r^m + o(r^m)$ as $r(s) \to 0$.
Since $r\dr f = \la \nf, r\dr\ra = pf$, we must have $p=m$ and so
$f(\nu) = -\ba \neq0)$.

The statement on $\xi_f$ holds since $\la \Hs (f)(\nu)\cdot \nu,\nu\ra = p(p-1)f(\nu)$.
\end{proof}

The following result should also be already known. Nevertheless, we state it and provide a proof.
\begin{proposition}\label{Proposition:hessian-linear-part}
Let $f= f_p + f_{p+1} + \ldots$, be the sum of the homogeneous parts of the real analytic function $f$.
Let $\nu$ be the limit of the secant of a trajectory $\absgamma$ of the gradient $\nf$ with
$\omega (\absgamma) =\{O\}$, and such that $f_p (\nu) = -\ba \neq 0$.
Then the divided gradient vector field $\xi_f$ is resolved at $\nu$, that is $\xi_f$ has a non-nilpotent
linear part at $\nu$. Moreover, the eigen-values of the linear part of $\xi_f$ are entirely
determined by $\Hs (f_p)(\nu)$ (which is also $\lim_{x/|x|\to\nu} |x|^{-p+2} \Hs (f) (x)$).
\end{proposition}
\begin{proof}
Let us look at the dominant part of the gradient vector field nearby the oriented direction $\nu$.
First we assume that $\nu = \partial_{x_1}$. Second if $H := \Hs (f_p)(\nu)$, and since $\nu$
is an eigen-vector of $H$, we take orthogonal coordinates $x=(x_1,y)$ at so that $H = [h_i]$
is diagonal, where $h_i$ is the $i$-th diagonal coefficient and $h_1 = -p(p-1)\ba \neq 0$. \\
Since $u_1 = \sqrt{1-\sum v_i^2}$, where $v_i = u_i$ for $i\geq 2$, we find that
$$\la \Hs (f_p)(u)\cdot u,u\ra = h_1 + \sum_{i\geq 2} (h_i - h_1)v_i^2.$$
Thus, in a neighborhood of
$\nu$ the divided gradient differential equation becomes
\begin{equation}
\left \{
\begin{array}{rcl}
\dot{r} & = &  h_1r + O(|(v,r)|^2) \\
\dot{v}_i & = & (h_i-h_1)v_i + O(|v|^3) +O(r),
\end{array}
\right.
\end{equation}
%
%
%
As announced we find that the eigen-values of the linear part of $\xi_f$ are easily deduced
from that of $\Hs (f_p) (\nu)$.
\end{proof}

The following is a straightforward consequence of Proposition \ref{Proposition:hessian-linear-part}.
\begin{corollary}\label{Corollary:hessian-eigen-almost-homogeneous}
Under the hypothesis of Proposition \ref{Proposition:hessian-linear-part}, we find
\begin{center}
\vspace{4pt}
$\displaystyle{\lim_{\absgamma\ni x \to O}} H (f)(x) = C \, \Hs (f_p) (\nu)$,
\vspace{4pt}
\end{center}
for a positive constant $C$.
\end{corollary}

\medskip
The second point we would like to put in light is the following: 

\smallskip\noindent
2-1) Assume that $f=f_p$ is a homogeneous polynomial of degree $p\geq 3$, let $\Sigma$ be the set of points
where the gradient $\nf$ is an eigen-vector of $\Hs (f)$, that is
\begin{center}
\vspace{4pt}
$\Sigma =\{x: \ud f \wedge \ud |\nf| = 0\} \subset \{x: \ud f \wedge \ud |\nf|^2 = 0\}$.
\vspace{4pt}
\end{center}
It is a closed semi-algebraic cone over $O$.\\
In particular each open ray $\R_{>0} u$, for $u \in (\Sigma \cap \Sr^{n-1})$ with $f(u)<0$,
is a gradient trajectory with $\omega$-limit point $O$.\\
For a positive number $r$, we define $\Sigma_r = \Sigma \cap f^{-1}(r)$. It has finitely many connected
components $\Sigma_r^1\ldots,\Sigma_r^q$ on each of which $|\nf|$ takes a constant values $pc_i (r)$.
Note that $\cup_{r>0} \Sigma_r^i$ is a positive open cone $C_i^*:=\R_{>0}\Sigma^i$ , where $\Sigma^i$ is a
closed semi-algebraic subset of $\Sr^{n-1}$. At any $x \in C_i^*$, using the spherical coordinates and 
writing $x=ru$, we get
\begin{center}
\vspace{4pt}
$|\nf (u)| = pr^{1-p}c_i (r) = pc_i$.
\vspace{4pt}
\end{center}
Up to a re-numeration we have $0<c_1\leq\ldots\leq c_q$.
At any $x =ru \in C_i^*$, we have $\Hs (f)(x)\cdot \nuf (x) = \lambda (x)\nu(x)$. We deduce
that $\lambda (x) = r^{p-2} \lambda (u)$, and so the mapping $\Sigma^i \ni u \to \lambda (u)$ 
is semi-algebraic.

If $\nu$ is a limit of secants at $O$ of a gradient trajectory $\absgamma$, then
we know that $\Hs (f)(\nu)\cdot\nu = \lambda\nu$ for $\lambda \neq 0$. So $\absgamma = \R_{>0} \nu$
and is contained in $\Sigma$. Since $\absgamma \subset C_i^*$, we get that $\nf \circ \gamma (r) =
- p r^{p-1} c_i \nu$, since we have agreed that $f\circ\gamma <0$. So we find that
\begin{center}
\vspace{4pt}
$\lambda (\nu) = p(p-1)c_i$ and so $f (\nu) = - c_i$.
\vspace{4pt}
\end{center}
This very special example suggests the existence of links between gradient trajectories, Hessian matrices 
and their eigen-values and limits of extrema of $|\nf|$ restricted to the foliation by the level of $f$.

\medskip\noindent
2-2)
Let $f$ be a real analytic Morse function, vanishing at the origin. We write $f = f_2 + O(|x|^3)$.
The divided gradient differential equation ensures that
the limit $\nu$ of the secants of any trajectory $\absgamma$ falling onto $O$ is an eigen-direction of
the Hessian matrix $\Hs (f_2)$, which is invertible at $O$.
Up to a an orthonormal change of coordinates we write

\begin{center}
\vspace{4pt}
$f (x) = -\sum_{i=1}^s \lambda_i x_i^2 + \sum_{j=1}^u\mu_jx_j^2+ O(|x|^3)$,
\vspace{4pt}
\end{center}
 where $s,u$ are non negative integers and $s+u =n$ and, when $s>0$ then $0<\lambda_1\leq\ldots\leq\lambda_s$ 
and, when $u>0$ then $0<\mu_1\leq\ldots\leq\mu_u$.
From Corollary \ref{Corollary:hessian-eigen-almost-homogeneous}, along a trajectory $\absgamma$
we get $r^{-1}\gamma \to \nu$ and that $\nu$ lies in the eigen-space corresponding to the eigen-value
$-2\lambda$ for $\lambda = \lambda_i$. Then we observe that along such a trajectory
$r^{-1}|\nf| \to 2\lambda$.

Let $\beta_1\leq\ldots\leq\beta_n$ be the ordered set of eigen-values counted with multiplicities.
We also check that for $\ve$ positive and small enough, the function $|\nf|_{| f^{-1}(-\ve)}$ has local
extrema $c_1(-\ve)\leq\ldots \leq c_q(-\ve)$, for $q\geq 1$, which are continuous subanalytic functions
of $\ve$. We also check that each such function $c(-\ve)\ve^{-1/2}$ tends to $2\beta$ as $\ve$ tends
to $0$ for $\beta = \beta_i$ for some $i$. 
\\
In other words the eigen-values $2\lambda_i$ and $2\mu_j$ are strongly connected with the limits
at $O$ of the local extrema of the restriction of $|\nf|$ to the (connected components of the) levels
of the function $f$.

%
%
%
%
%
%
%
%
%
%
%
%
%
%
%
%
%
%
%
%
%
%
%
%
%
%
%
\section{Main result}\label{Section:main-result}
We use the notation introduced in Section \ref{Section:setting}.

\medskip
\noindent
{\bf Non-oscillating hypothesis:} We assume that $\absgamma$ does not oscillate at its end point $O$.

\medskip
We recall that plane gradient trajectories do not oscillate at their limit point, and that
in higher dimension non-oscillating gradient trajectories exist as well \cite{Mo,Sa,FS,GS}.

\medskip
We begin with an elementary but key result about the  monotonicity of functions along a non oscillating
trajectory.
\begin{lemma} \label{Lemma:monotonic}
Let $\psi: U \to \R$ be a continuous semi-analytic function, analytic and non zero outside $\crit (f)$.
The function $\psi|_{\absgamma}$ has a finite limit at $O$, and up to shrinking $U$, the function
$\psi|_{\gamma}$ is either constant or strictly monotonic on $\absgamma\cap U$.
\end{lemma}
\begin{proof}
Assume $\psi|_{\absgamma}$ is not constant.

Assume there exists a sequence $(s_j)_j$
converging to $l(\absgamma)$, the length of $\absgamma$ such that $\psi^\prime$ vanishes and changes sign at
each $s_j$. Since
\begin{center}
\vspace{4pt}
$\psi^\prime (s) = \la \nabla \psi,\nuf\ra \circ \gamma (s)$,
\vspace{4pt}
\end{center}
the subset $\absgamma \cap \{ \la \nabla \psi,\nuf\ra = 0\}$ has thus infinitely many connected components.
Since $\clos(\{ \la \nabla \psi,\nuf\ra = 0\})$ is a closed semi-analytic set and
contains $O$, thus this fact contradicts the non-oscillation of $\absgamma$.
\end{proof}

The function $s \to r = |\gamma (s)|$ is strictly decreasing to $0$ as $s\to l(\Gamma)$, following
Lemma \ref{Lemma:monotonic}. Moreover we have
\begin{center}
\vspace{4pt}
$\dot{r} (s) = \displaystyle{\frac{|\nf|}{\drf} \circ \gamma (s)}$ and $\dot{r}(s) \to -1$ as $s \to
l(\Gamma)$.
\vspace{4pt}
\end{center}
So we reparameterize the integral curve $\gamma$ with $r$:
\begin{center}
\vspace{4pt}
$\gamma^\prime (r) = \displaystyle{\frac{\ud}{\ud r} (\gamma(r)) = \left (\frac{|\nf|}{\drf} \;  \nuf \right )
\circ \gamma (r)} =\xi \circ \gamma (r)$,
\vspace{4pt}
\end{center}
where $\xi$ is the semi-analytic vector field appearing in the above expression.

Another simple but key result is the following
\begin{lemma}\label{Lemma:main2}
Let $a,b$ be  continuous semi-analytic functions with $\{b=0\} \subset \{\drf =0\}$.
Let $v$ be the function defined as $v(r) := (a/b) \circ \gamma (r)$. Assume that $v$ is bounded.
Thus $r v^\prime (r) \to 0$ as $r\to 0$.
\end{lemma}
\begin{proof}
First note that by the nature of  $\gamma (r)$ we will still have $v^\prime = \bar{a}/b^2$ having the 
same form as $v$.
From the proof of Lemma \ref{Lemma:monotonic}, we know that close enough to $0$, the function $v$ is 
strictly monotonic
or constant. If $v^\prime \to \infty$, then $rv^\prime$ is still strictly monotonic close enough to $0$.
So if $rv^\prime$ does not tend to $0$, integrating $v^\prime$ would contradict the fact that $v$ is bounded.
\end{proof}

A first interesting consequence of Lemmas \ref{Lemma:monotonic} and \ref{Lemma:main2}.
\begin{corollary}\label{Corollary:limit-tangent}
$\lim_{r\to 0}\nu_f \circ  \gamma (r) = -\nu$.
\end{corollary}
\begin{proof}
First we write $\gamma (r) =  r(\nu+\chi (r))$, where $\chi (r) = (\dr -\nu)\circ \gamma (r)$.
Since the vector field $\dr$ is semi-analytic, applying \ref{Lemma:monotonic} we deduce that such
that $|\chi(r)|\to 0$ as $r\to 0$. Thus we deduce $\gamma^\prime (r) =  s^\prime(r)\nu+\beta(r)$
where $\beta (r) = \chi (r)+ r\chi^\prime (r)$.
From Lemma \ref{Lemma:monotonic} and Lemma \ref{Lemma:main2}, we get $\beta (r) = o(1)$.
And to conclude we recall that $s^\prime (r) \to -1$ as $r\to 1$,
\end{proof}

\medskip
We recall that $H(f)$ is \em "the normalized Hessian matrix" \em of $f$, and so the mapping $x\to H(f)(x)$
 continuous and semi-analytic. In particular thanks to Lemma \ref{Lemma:monotonic} we deduce the following
\begin{lemma}\label{Lemma:limit-hessian}
There exists $\Cal{H} \in \Symn^*$ such that
\begin{center}
\vspace{4pt}
$\lim_{\absgamma\ni x \to O}H(f)(x) = \Cal{H}$.
\vspace{4pt}
\end{center}
\end{lemma}
Note that this result can also be deduced from the results of \cite{KPau}.

\smallskip\noindent
This fact satisfied, the main result of this note is:
\begin{theorem}\label{Theorem:main}
$$\lim_{\absgamma\ni x \to O} \frac{\Hs(f)(x)\cdot \nu_f (x)}{|\Hs(f)(x)\cdot \nu_f (x)|} = 
\lim_{\absgamma\ni x \to O} \frac{H(f)(x)\cdot \nu_f (x)}{|H(f)(x)\cdot \nu_f (x)|} = 
-\nu.$$
Thus the vector $\nu$ is an eigen-vector of the matrix $\Cal{H}$.
\end{theorem}
The proof will consist of several elementary results.

\begin{proof}
In the proof of Corollary \ref{Corollary:limit-tangent}, we had written 
$\gamma (r) = r(\nu+\chi(r))$, so that $\gamma^\prime =  s^\prime(r)\nu + \beta$ with 
$\beta = \chi + r\chi^\prime$. Since
\begin{center}
$\chi = \dr -\nu\;\;$ and $\;\;\chi^\prime (r) = \displaystyle{\frac{1}{r\dr f}\np f \circ \gamma (r)}$,
\end{center}
we can write $\beta^\prime (r) = (a/b) \circ \gamma (r)$ for a semi-analytic function $v=a/b$ with
$b^{-1}(0) =\dr f^{-1}(0)$. Applying again Lemma \ref{Lemma:main2}, we deduce
\vspace{4pt}
\begin{equation}\label{Equation:eq-1}
\gamma^{\prime\prime} (r) = o(r^{-1}).
\vspace{4pt}
\end{equation}
Note that
\begin{center}
$\gamma^{\prime\prime} (r) = v^\prime (r)\nu_f(\gamma(r))+ v^2(r)( \ud \nuf \cdot \nuf)
\circ \gamma (r)\;$ with $\;v(r) = \displaystyle{\frac{|\nf|}{\drf}}\circ \gamma (r)$.
\end{center}
We recall that $v(r) \to -1$ as $r\to 0$ and so from Lemma \ref{Lemma:main2}, $r v^\prime (r) \to 0$ as
$r\to 0$. We deduce from (\ref{Equation:eq-1})
\vspace{4pt}\begin{equation}\label{Equation:eq-2}
(\ud \nuf \cdot \nuf) \circ \gamma (r)= o(r^{-1}).
\vspace{4pt}
\end{equation}

Since
\begin{center}
$\ud \nuf \cdot \nuf  =
\displaystyle{\frac{1}{|\nf|}}[\Hs f \cdot \nuf -\langle \Hs f\cdot\nuf,\nuf\rangle\nuf]$.
\end{center}
Equation (\ref{Equation:eq-2}) now reads
\vspace{4pt}
\begin{equation}\label{Equation:eq-3}
r[\Hs f \cdot \nuf -\langle \Hs f\cdot\nuf,\nuf\rangle\nuf] \circ \gamma (r) = o (|\nf|).
\vspace{4pt}
\end{equation}
The classical Bochnak-\L ojasiewicz inequality applied to $f$ and then to $|\nf|$ nearby the origin $O$ gives
\begin{center}
\vspace{4pt}
 $|x|\cdot | \Hs f\cdot\nf| \geq C|\nf|^2$ for a positive constant $C$.
\vspace{4pt}
\end{center}
or equivalently
\begin{center}
\vspace{4pt}
 $|x|\cdot | \Hs f\cdot\nu_f| \geq C|\nf|$ for a positive constant $C$.
\vspace{4pt}
\end{center}
Since $r|\Hs f| \cdot \nuf \geq C|\nf|$ the only possibility for Equation (\ref{Equation:eq-3})
to hold is that along $\absgamma$ we get
\vspace{4pt}
\begin{equation}\label{Equation:eq-4}
\Hs f \cdot \nuf = \langle \Hs f\cdot\nuf,\nuf\rangle\nuf + o(| \langle \Hs f\cdot\nuf,\nuf\rangle|).
\vspace{4pt}
\end{equation}
Dividing both sides by $|\Hs (f)\nu_f|$ and using Corollary \ref{Corollary:limit-tangent} 
provides that the limit is either $\nu$ or $-\nu$.
Along $\absgamma$, the estimates $f = -\ba r^m\dr + o(r^m)$ from Proposition \ref{Proposition:exponent-vca} 
will give that it is $-\nu$.  
\end{proof}
%
%
%
%
%
%
%
%
%
%
%
%
%
%
%
%
%
%
%
%
%
%
%
%
%
%
%
%
%
%
%
%
%
%
\section{A few consequences}\label{Section:AFUC}
For a given non-oscillating trajectory $\absgamma$ with $\omega$-limit point $O$,
we deduce from Proposition \ref{Proposition:exponent-vca}, Lemma \ref{Lemma:monotonic} and 
Lemma \ref{Lemma:main2} that along $\absgamma$ we have
\vspace{4pt}
\begin{equation}\label{Equation:gradient-radial}
\nabla f = -m \ba r^{m-1}\dr + o(r^{m -1}).
\vspace{4pt}
\end{equation}

The first un-expected result is the following one.

\begin{proposition}\label{Proposition:vca-gradient}
Let $\absgamma$ be a non oscillating gradient trajectory of the function $f$ with $\omega$-limit point
the origin $O$. Let $m$ and $\ba$ be as above. \\
The real number $m\cdot\ba$ is an asymptotic critical value of the function $r^{1-m} |\nabla f|$.
\end{proposition}
\begin{proof}
Let us denotes by $R$ the function $r\to r^{1-m} |\nabla f(\gamma(r))|$. Since $R \to m\ba$ as $r\to 0$,
Lemma \ref{Lemma:main2} implies  $rR^\prime \to 0$ as $r\to 0$.

Since
\begin{center}
$
\displaystyle{ R^\prime (r) = \frac{1}{r^m} \frac{\ud s}{\ud r}\langle r
\Hs f\cdot \nu_f - (m-1) |\nf|\dr ,  \nu_f \rangle}
$
\end{center}
and $\frac{\ud s}{\ud r} \to -1$ we get
\begin{center}
\vspace{4pt}
$\langle r \Hs f\cdot \nu_f - (m-1) |\nf|\dr ,  \nu_f \rangle = o(r^{m-1})$.
\vspace{4pt}
\end{center}
Along $\absgamma$ we know that $\nf \simeq \ba r^{m-1} \dr$ and so by Bochnak-\L ojasiewicz
$r|\Hs f\cdot\nu_f| \geq C r^{m-1}$.
From Theorem \ref{Theorem:main} we also have $|\langle \Hs f\cdot\nu_f, \nu_f \rangle| \simeq |\Hs f \cdot
\nu_f|$, and so we necessarily deduce that
\begin{center}
\vspace{4pt}
$r \Hs(f) \cdot \nu_f = (m-1) |\nf| \dr + o(r^{m-1})$.
\vspace{4pt}
\end{center}
So we just deduce that along $\absgamma$ the gradient vector field $\nabla |\nabla f|$ is asymptotically
radial. Thus the result is proved.
\end{proof}

There are two consequences summarized in the next statements.
The first one is already in the proof of Proposition \ref{Proposition:vca-gradient}:
\begin{corollary}\label{Corollary:consequence-1}
Along $\absgamma$, we find
\begin{center}
\vspace{4pt}
\begin{tabular}{rcl}
\vspace{4pt}
$\Hs(f) \cdot \nu_f$ &  = & $ m (m-1) \ba r^{m-2} \nu_f + o(r^{m-2})$ \\
\vspace{4pt}
& = & $  -m (m-1) \ba r^{m-2} \dr + o(r^{m-2})$ \\
\vspace{4pt}
& = & $  m (m-1) \ba r^{m-2} \nu + o(r^{m-2})$.
\vspace{4pt}
\end{tabular}
\end{center}
\end{corollary}
Let $g$ be a homogeneous polynomial of degree $p\geq 2$. We recall that
\vspace{4pt}
\begin{equation}\label{Equation:hessian-homogeneous}
p (p-1) g= \langle \Hs (g) \cdot r\dr, r\dr\rangle,
\vspace{4pt}
\end{equation}

For Corollary \ref{Corollary:consequence-1} and Lemmas \ref{Lemma:monotonic} and \ref{Lemma:main2}
we see that along $\absgamma$, we obtain an estimate similar to Equation 
(\ref{Equation:hessian-homogeneous}):
\begin{center}
\vspace{4pt}
$m (m-1) f = \langle \Hs f \cdot r\dr, r\dr\rangle + o(r^m)$.
\vspace{4pt}
\end{center}

\medskip
The second consequence is that the asymptotic critical value $-\ba$ may be interpreted as an eigen-value
under some additional hypotheses:
\begin{corollary}\label{Corollary:consequence-2}
Assume that along $\absgamma$, $r^2 h (\gamma(r)) \sim r^m \sim |f|$.
Then along $\absgamma$ the matrix $\frac{1}{m(m-1) r^{m -2}} \Hs (f)$ has a limit at $O$ in
$\Symn^*$ and the eigen-value in the oriented direction $\R_+ \nu$ is exactly
$- \ba$.
\end{corollary}
%
%
%
%
%
%
%
%
%
%
%
%
%
%
%
%
%
%
%
%
%
%
%
%
%
%
%
%
%
%
%
%
%
%
%
%
%
\section{At infinity}\label{Section:infinity}
Even though the results of \cite{Gr1} were obtained in the semi-algebraic context, they still hold true
for $C^2$ globally subanalytic function.\\
Let $f:\Rn\to \R$ be a globally subanalytic real analytic function.
We consider again the gradient differential equation:
\begin{center}
\vspace{4pt}
$\dot{\gamma}(s) = \nuf (\gamma (s))$ with $\gamma (0) = x_0$ and $f(x_0) < 0$,
\vspace{4pt}
\end{center}
But we are now interested in trajectories that leave any compact subset, that is, for $s\to \gamma(s)$ a
parameterized solution, given any compact subset $K$ of $\Rn$, there exists $s_K$, such that $\gamma(s) 
\notin K$ for $s>s_K$, or equivalently, $|\gamma (s)|\to +\infty$ as $s \to \infty$.

\smallskip\noindent
Assume that $\absgamma$ is a non-bounded gradient trajectory, then $\absgamma$ leaves any compact.
We still suppose that $f$ is
increasing along $\absgamma$. We may assume that the $\lim_{x\in\absgamma,|x|\to +\infty} f(x) = 0$.

\smallskip\noindent
The main result of \cite{Gr1} is that the gradient conjecture is also true at infinity:
\begin{theorem}[\cite{Gr1}]\label{Theorem:gradient-conjecture-infinity}
There exists $\nu \in \Sr^{n-1}$ such that
$\displaystyle{\lim_{x\in\absgamma, |x|\to +\infty}\frac{x}{|x|} =  \nu}$.
\end{theorem}

To prove this result, some estimates along the trajectory were necessary. We recall one that is relevant for
this section
\begin{proposition}[\cite{Gr1}]\label{Proposition:exponent-vca-infinity}
There exists a positive rational number $m$ and a positive real number $\ba$ such that along $\absgamma$
\begin{center}
\vspace{4pt}
$ f = -\ba r^{-m} + o(r^{-m})$ and $\dr f =   m\ba r^{-1-m} + o(r^{-1-m})$
\vspace{4pt}
\end{center}
\end{proposition}

We need the following definition
\begin{definition}
The trajectory $\absgamma$ is {\rm non-oscillating at infinity} if for any globally subanalytic
non-bounded analytic hypersurface $H$, there exists a radius $R$ such that $\absgamma\cap H \cap
(\Rn\setminus \bB(O,R))$ has finitely many
connected components.
\end{definition}

The main result in this context, using the notations introduced in the previous sections, is
\begin{theorem}\label{Theorem:main-infinity} Assume that $\absgamma$ is non-oscillating at infinity.
Then

\vspace{4pt}
\noindent
1) $\lim_{x\in\absgamma, |x|\to +\infty} H(f)(x) = \Cal{H} \in  \Symn^*$

\vspace{4pt}
\noindent
2) The vector $\nu$ is an eigen-vector of the matrix $\Cal{H}$.
\end{theorem}

The proof will work along the same lines with a few modifications that we will describe.
The first one is that there exists a Bochnak-\L ojasiewicz inequality at infinity:
\begin{proposition}[\cite{DDVG1}]\label{Proposition:BL-infinity}
If $|x| \gg 1$ and $|f(x)| \ll 1$ then
\begin{center}
\vspace{4pt}
$|x|\cdot |\nf| \geq C|f|$,
\vspace{4pt}
\end{center}
for a positive constant $C$.
\end{proposition}
Then using it for the function $|\nf|$ we find again that for $|x| \ll 1$ and $|\nf| \ll 1$, there is a
positive constant $C$.
\vspace{4pt}
\begin{equation}\label{Equation:BL-infinity}
|x|\cdot |\Hs f\cdot \nabla f| \geq C|\nabla f|^2.
\vspace{4pt}
\end{equation}
Next we need the analogues of Lemmas \ref{Lemma:monotonic} and \ref{Lemma:main2}.

\begin{lemma}\label{Lemma:main-infinity-1}
Let $g$ be a globally continuous subanalytic function outside of a compact subset of $\Rn$.
Assume that $\absgamma$ is non oscillating at infinity.
Let $v := g\circ \gamma$. As a germ at infinity, if $v$ is not constant then it is strictly monotonic.
\end{lemma}
\begin{proof}
If not it is contradiction to the non-oscillating at infinity of $\absgamma$.
\end{proof}

Since the function $1/r$ is semi-algebraic, we can then reparameterize
the solution $\gamma$ by $r = |\gamma (s)|$.
\begin{lemma}\label{Lemma:main-infinity-2}
Let $a,b$ be  globally continuous subanalytic function outside of a compact subset of $\Rn$, with
$\{b=0\} \subset \{\drf =0\}$.
Assume that $\absgamma$ is non oscillating at infinity.
Let $v (r):= (a/b)\circ \gamma (r)$ and assume that $v (r) \to c\in\R$ as $r\to +\infty$. Then the function
$r\mapsto r^{-1} v^\prime(r)$ tends to $0$ as $r\to +\infty$.
\end{lemma}
\begin{proof}
Similar to that of  Lemmas \ref{Lemma:main2}.
\end{proof}

From these two Lemmas, and adapting the proof Corollary \ref{Corollary:limit-tangent}, it is  easy to prove
the following
\begin{corollary}\label{Corollary:limit-tangent-infinity}
Assuming that $\absgamma$ is non oscillating at infinity, we get
\begin{center}
\vspace{4pt}
$\lim_{r\to 0}\nu_f \circ  \gamma (r) = \nu$.
\end{center}
\end{corollary}

The proof of Theorem \ref{Theorem:main-infinity} goes along the same lines as those of the proof of
Theorem \ref{Theorem:main-infinity} using Lemmas \ref{Lemma:main-infinity-1} and \ref{Lemma:main-infinity-2}.
Corollary \ref{Corollary:limit-tangent-infinity} in a similar way.

\medskip
Along a non-oscillating trajectory $\absgamma$ going to infinity, we find estimates  similar to Equation
(\ref{Equation:gradient-radial}), Proposition \ref{Proposition:vca-gradient}, Corollary
\ref{Corollary:consequence-1} and Corollary \ref{Corollary:consequence-2}, namely
\vspace{4pt}
\begin{equation}\label{Equation:gradient-radial}
\nabla f = m \ba r^{-m-1}\dr + o(r^{-m-1}).
\vspace{4pt}
\end{equation}
We also find
\begin{proposition}\label{Proposition:vca-gradient-infinity}
Let $\absgamma$ be a gradient trajectory of the function $f$ non-oscillating at infinity.
Let $m$ and $\ba$ be as above. \\
The real number $m\cdot\ba$ is an asymptotic critical value of the function $r^{1-m} |\nabla f|$.
\end{proposition}
\begin{proof}


It works as in the proof of Proposition \ref{Proposition:vca-gradient}
\end{proof}

We eventually  deduce the following estimates
\begin{center}
\vspace{4pt}
\begin{tabular}{rcl}
\vspace{4pt}
$\Hs(f) \cdot \nu_f$ &  = & $- m (m-1) \ba r^{-m-2} \nu_f + o(r^{-m-2})$ \\
\vspace{4pt}
& = & $- m (m-1) \ba r^{-m-2} \dr + o(r^{-m-2})$ \\
\vspace{4pt}
& = & $- m (m-1) \ba r^{-m-2} \nu + o(r^{-m-2})$ \\
\vspace{4pt}
$m (m-1) f$ & = & $\langle \Hs f \cdot r\dr, r\dr\rangle + o(r^{-m})$.
\end{tabular}
\vspace{4pt}
\end{center}
%
%
%
%
%
%
%
%
%
%
%
%
%
%
%
%
%
%
%
%
%
%
%
%
%
%
%
%
%
%
%
%
%
%
%
%
%
%
%
%
%
%
%
%
\section{Riemannian case}\label{Section:riemannian}

Before getting in what really interest us, we will recall a few basic formulae of 
Differential Geometry.

Let $(M,\bg)$ be a smooth Riemannian manifold of dimension $n$. Let $TM$ be its tangent 
bundle and $\Scr{X}(M)$ be the $C^\infty(M)$-module of the smooth vector fields on $M$.
Let $\la , \ra_\bg$ be the scalar product coming with $\bg$ and $|-|_\bg$ the associated 
norm. Let $D$ be the canonical connection associated with $\bg$. Given $X,Y,Z \in\Scr{X}(M)$, 
the connection $D$ satisfies the following equations:
\begin{equation}
[X,Y] = D_XY - D_YX= \ud Y\cdot X - \ud X\cdot Y,
\end{equation}
\begin{equation}
X\cdot\la Y,Z\ra_\bg = \la D_XY,Z\ra_\bg + \la Y,D_XZ\ra_\bg,
\end{equation}
\begin{equation}
D_XY = \ud Y\cdot X + \Gamma (X,Y),
\vspace{4pt}
\end{equation}
Where $\Gamma:\Scr{X}(M)\times\Scr{X}(M)\to\Scr{X}(M)$ is smooth, $C^\infty$-bilinear and symmetric.
In local coordinates $\Gamma$ is defined as usual with the mean of Christoffels symbols.

\medskip
Given a smooth function $u$, we recall that the gradient vector field $\ngu$ is defined as the ``dual''
of $Du$ for the pairing $\la,\ra_\bg$, that is for any smooth vector field $X$ on $M$,
\vspace{4pt}
\begin{equation}
(Du)(X) = D_Xu = X\cdot u = \la\ngu,X\ra_\bg.
\end{equation}

The notion of Hessian of the function $u$ as $\ud^2 u$ is no longer accurate in this Riemannian setting.
Since the connection $D$ extends to tensors, the corresponding notion of Hessian in this setting comes from
the following tensor
\begin{center}
\vspace{4pt}
$D^2u = D(Du)=(D\ud u):\Scr{X}(M)\times\Scr{X}(M)\to C^\infty(M)$

\end{center}
defined as
\begin{center}
$(D\ud u)(X,Y) = (D_X \ud u)(Y)$.
\vspace{4pt}
\end{center}
 It is smooth, $C^\infty$-bilinear and symmetric as well.
Using the formulae above, we also check that
\vspace{4pt}
\begin{equation}
(D\ud u)(X,Y) =  \la D_X\ngu,Y\ra_\bg = \la D_Y\ngu,X\ra_\bg,
\vspace{4pt}
\end{equation}
and deduce that
\vspace{4pt}
\begin{equation}
\nabla_{\!\bg} |\ngu|_\bg^2 = 2D_\ngu \ngu.
\end{equation}
We also look at $D\ud u$ as a $C^\infty(M)$-linear map $\Scr{X}(M)\to\Scr{X}(M)$, defined, with an obvious
abuse of notation, as $(D\ud u)(X) = D_X\ngu$.
At each point $x$ of $M$ the tensor $D\ud u$ induces a bilinear symmetric form
$\Hsg (u)(x): =(D\ud u)_x =T_xM\times T_xM\to \R$.
We will call it the Hessian of $u$ at $x$ and obviously $x \to \Hsg (u)(x)$ is smooth. And we will also
look at it as a linear endomorphism of $T_xM$.

\bigskip
Now we come to our topic.

Assume that $M=U$ is an open subset of $\Rn$ containing the origin $O$.

Assume that the Riemannian metric $\bg$ is analytic.

Let $f:U \to \R$ be a real analytic function such that $O \in \crit (f)$ the critical locus of $f$.
We also assume $f(O) = 0$.

Let us denote $\nugf$ the unitary gradient vector field defined on $U\setminus \crit (f)$.

Let us consider the following gradient differential equation:
\begin{center}
\vspace{4pt}
$\dot{\gamma}(s) = \nugf (\gamma (s))$ with $\gamma (0) = x_0$ and $f(x_0) < 0$.
\vspace{4pt}
\end{center}

We have picked some orthonormal coordinates at $O$ so that $\bg (O)$ is the Euclidean metric.

If $r$ is the Euclidean distance function to the origin $O$ and $r^\bg$ the distance to the origin for
the metric $\bg$, then, and $r^\bg/r \to 1$ as $x\to O$.
\begin{remark}
In a neighborhood of $O$, Bochnak-\L ojasiewicz inequalities holds,  namely
\begin{center}
\vspace{4pt}
$r^\bg | \ngf|_\bg \geq C |f|$ and so $r |\ngf|_\bg \geq C^\prime |f|$,
\vspace{4pt}
\end{center}
 for some positive constant $C,C^\prime$.
\end{remark}

Moreover the results of \cite{KMP} quoted in Section \ref{Section:setting} still hold true, sticking
with $r$ in the estimates.

Bochnak-\L ojasiewicz inequality for the function $|\ngf|$, then reads in a neighborhood of $O$
\begin{center}
\vspace{4pt}
$r |D_\nugf\ngf|_\bg \geq C |\ngf|_\bg$ for some positive constant $C$.
\vspace{4pt}
\end{center}
The operator $\Hsg (f)$ is now  analytic and in the present coordinates, when we look at it as linear mapping,
this gives us (symmetric) analytic entries $h_{i,j}$. We recall that  the semi-analytic function $h$ is
$h = (\sum_{i,j}h_{i,j}^2)^{1/2}$ and we define $H^\bg(f)  = h^{-1}\Hsg (f)$, so that at each point
$x\notin\crit (f)$ we deduce $H^\bg (f) (x)$ is non zero.

\medskip
Now we can state the main result of this section in this general case.
\begin{theorem}\label{Theorem:main-riemmannian} Assume that $\absgamma$ is a non-oscillating trajectory
of the gradient field associated with $f$ respectively to the analytic Riemannian metric $\bg$,
with $\omega (\absgamma) = \{O\}$.
Then

\vspace{4pt}
1) Along $\absgamma$, $\lim_{x\in\absgamma, |x|\to +\infty} H^\bg (f) (x) = \Cal{H}^\bg \in  {\rm Sym}(\Rn)
\setminus \bo$

\vspace{4pt}
2) $$\lim_{\absgamma\ni x \to O} \frac{\Hsg(f)(x)\cdot \nu_f (x)}{|\Hsg(f)(x)\cdot \nu_f (x)|} = 
\lim_{\absgamma\ni x \to O} \frac{H^\bg(f)(x)\cdot \nu_f (x)}{|H^\bg(f)(x)\cdot \nu_f (x)|} =-\nu.$$
The vector $\nu$ is an eigen-vector of the matrix $\Cal{H}^\bg$.
\end{theorem}
\begin{proof}
The first point is elementary.
The proof of the second point goes exactly along the same steps as in the Euclidean case dealt with in
Section \ref{Section:main-result}. We nevertheless has to check the form of $\ddot{\gamma}$ in order
to recycle the previous proof.
\begin{center}
\vspace{4pt}
$\displaystyle{\ddot{\gamma}(s) = \ud\nugf\cdot\nugf}$.
\vspace{4pt}
\end{center}
Thus we find
\begin{center}
\vspace{4pt}
$\displaystyle{\ud \nugf \cdot \nugf = \frac{1}{|\ngf|_\bg}[D_\nugf\ngf - \Gamma (\nugf,\ngf) -
\la D_\nugf \ngf,\nugf\ra_\bg \nugf]}$.
\vspace{4pt}
\end{center}
and we observe that $\Gamma (\nugf,\ngf) = O(|\ngf|)$. Thus applying Bochnak-\L ojasiewicz for $D_\nugf \ngf$
will allow to conclude again.
\end{proof}

\begin{remark}
1) It is possible to find estimates as in Section \ref{Section:AFUC} but it is not clear what information
some of them will really carry. \\
2) In this context we will not speak of what is happening at infinity since the paper \cite{Gr1} only dealt
with the Euclidean metric. Moreover we would need to know some facts about the behavior of the metric
$\bg$ at infinity and it is clear what should be required to make things work.
\end{remark}
%
%
%
%
%
%
%
%
%
%
%
%
%
%
%
%
%
%
%
%
%
%
%
\section{Remarks and comments}\label{Section:remarks-and-comments}
\noindent
{\bf First.} What we have proved is not so much attached to a function than to a foliation,
namely the foliation by the (germ at $O$ of the) connected components of the levels of the given function.
In particular, there is a unique "matrix" in the space of directions $\bsym (n)$ which is attached, at
the limit point, to a given non-oscillating gradient trajectory. What this object encodes is still very
unclear. We have clearly suggested that it is linked with the normal form of the divided gradient vector
field, but also with the lines of ridges and valleys of the given foliation.  \\
At this stage, it is not guaranteed that what we just suggested carries the relevant phenomena(on) to
understand better what we have uncovered here.

\medskip\noindent
{\bf Second.}
Assume we are given an o-minimal structure $\Cal{M}$ expanding the real numbers in which the semi-analytic
sets are definable.
We check that for $C^2$ plane function definable in $\Cal{M}$ in a general $C^2$ definable
Riemannian setting our main results is still true. In the Euclidean setting and at infinity we need
globally definable to have it correct.\\
If the o-minimal structure $\Cal{M}$ is polynomially bounded, everything we proved in this note
is still true for any $C^2$ definable function (globally definable at infinity), even the estimates
The main reason for this is that Thom's Gradient conjecture is still true in these categories \cite{KPar}.
Note that in the plane, gradient trajectories are still definable in an o-minimal structure: The
Pfaffian closure of the initial o-minimal structure $\Cal{M}$.

\medskip\noindent
{\bf Third.} Proposition \ref{Proposition:hessian-linear-part} suggests that
the limit $\Cal{H}$ of our main result is linked with the reduction of the divided gradient vector field.
Corollary \ref{Corollary:consequence-2} suggests that at such a point $\nu$ the divided gradient vector field
may be simple, but we have no idea whether it is true or not and then how to prove/dis-prove this.
This is to be compared with the main result of \cite{CMR}, where is proved a resolution/reduction
result at the end point of a subanalytically non-oscillating trajectory of real analytic vector field in a
$3$-manifold.

\end{document}